\magnification=1200

\hsize=11.25cm    
\vsize=18cm     
\parindent=12pt   \parskip=5pt     

\hoffset=.5cm   
\voffset=.8cm   

\pretolerance=500 \tolerance=1000  \brokenpenalty=5000

\catcode`\@=11

\font\eightrm=cmr8         \font\eighti=cmmi8
\font\eightsy=cmsy8        \font\eightbf=cmbx8
\font\eighttt=cmtt8        \font\eightit=cmti8
\font\eightsl=cmsl8        \font\sixrm=cmr6
\font\sixi=cmmi6           \font\sixsy=cmsy6
\font\sixbf=cmbx6

\font\tengoth=eufm10 
\font\eightgoth=eufm8  
\font\sevengoth=eufm7      
\font\sixgoth=eufm6        \font\fivegoth=eufm5

\skewchar\eighti='177 \skewchar\sixi='177
\skewchar\eightsy='60 \skewchar\sixsy='60

\newfam\gothfam           \newfam\bboardfam

\def\tenpoint{
  \textfont0=\tenrm \scriptfont0=\sevenrm \scriptscriptfont0=\fiverm
  \def\rm{\fam\z@\tenrm}
  \textfont1=\teni  \scriptfont1=\seveni  \scriptscriptfont1=\fivei
  \def\oldstyle{\fam\@ne\teni}\let\old=\oldstyle
  \textfont2=\tensy \scriptfont2=\sevensy \scriptscriptfont2=\fivesy
  \textfont\gothfam=\tengoth \scriptfont\gothfam=\sevengoth
  \scriptscriptfont\gothfam=\fivegoth
  \def\goth{\fam\gothfam\tengoth}
  
  \textfont\itfam=\tenit
  \def\it{\fam\itfam\tenit}
  \textfont\slfam=\tensl
  \def\sl{\fam\slfam\tensl}
  \textfont\bffam=\tenbf \scriptfont\bffam=\sevenbf
  \scriptscriptfont\bffam=\fivebf
  \def\bf{\fam\bffam\tenbf}
  \textfont\ttfam=\tentt
  \def\tt{\fam\ttfam\tentt}
  \abovedisplayskip=12pt plus 3pt minus 9pt
  \belowdisplayskip=\abovedisplayskip
  \abovedisplayshortskip=0pt plus 3pt
  \belowdisplayshortskip=4pt plus 3pt 
  \smallskipamount=3pt plus 1pt minus 1pt
  \medskipamount=6pt plus 2pt minus 2pt
  \bigskipamount=12pt plus 4pt minus 4pt
  \normalbaselineskip=12pt
  \setbox\strutbox=\hbox{\vrule height8.5pt depth3.5pt width0pt}
  \let\bigf@nt=\tenrm       \let\smallf@nt=\sevenrm
  \normalbaselines\rm}

\def\eightpoint{
  \textfont0=\eightrm \scriptfont0=\sixrm \scriptscriptfont0=\fiverm
  \def\rm{\fam\z@\eightrm}
  \textfont1=\eighti  \scriptfont1=\sixi  \scriptscriptfont1=\fivei
  \def\oldstyle{\fam\@ne\eighti}\let\old=\oldstyle
  \textfont2=\eightsy \scriptfont2=\sixsy \scriptscriptfont2=\fivesy
  \textfont\gothfam=\eightgoth \scriptfont\gothfam=\sixgoth
  \scriptscriptfont\gothfam=\fivegoth
  \def\goth{\fam\gothfam\eightgoth}
  
  \textfont\itfam=\eightit
  \def\it{\fam\itfam\eightit}
  \textfont\slfam=\eightsl
  \def\sl{\fam\slfam\eightsl}
  \textfont\bffam=\eightbf \scriptfont\bffam=\sixbf
  \scriptscriptfont\bffam=\fivebf
  \def\bf{\fam\bffam\eightbf}
  \textfont\ttfam=\eighttt
  \def\tt{\fam\ttfam\eighttt}
  \abovedisplayskip=9pt plus 3pt minus 9pt
  \belowdisplayskip=\abovedisplayskip
  \abovedisplayshortskip=0pt plus 3pt
  \belowdisplayshortskip=3pt plus 3pt 
  \smallskipamount=2pt plus 1pt minus 1pt
  \medskipamount=4pt plus 2pt minus 1pt
  \bigskipamount=9pt plus 3pt minus 3pt
  \normalbaselineskip=9pt
  \setbox\strutbox=\hbox{\vrule height7pt depth2pt width0pt}
  \let\bigf@nt=\eightrm     \let\smallf@nt=\sixrm
  \normalbaselines\rm}

\tenpoint

\def\pc#1{\bigf@nt#1\smallf@nt}         \def\pd#1 {{\pc#1} }

\frenchspacing

\def\raggedbottom{\topskip 10pt plus 36pt\r@ggedbottomtrue}

\def\pointir{\unskip . --- \ignorespaces}

\def\Medbreak{\vskip-\lastskip\medbreak}

\long\def\th#1 #2\enonce#3\endth{
   \Medbreak\noindent
   {\pc#1} {#2\unskip}\pointir{\it #3}\smallskip}

\def\decale#1{\smallbreak\hskip 28pt\llap{#1}\kern 5pt}
\def\decaledecale#1{\smallbreak\hskip 34pt\llap{#1}\kern 5pt}
\def\puce{\smallbreak\hskip 6pt{$\scriptstyle\bullet$}\kern 5pt}

\def\eqalign#1{\null\,\vcenter{\openup\jot\m@th\ialign{
\strut\hfil$\displaystyle{##}$&$\displaystyle{{}##}$\hfil
&&\quad\strut\hfil$\displaystyle{##}$&$\displaystyle{{}##}$\hfil
\crcr#1\crcr}}\,}

\catcode`\@=12

\showboxbreadth=-1  \showboxdepth=-1

\newcount\numerodesection \numerodesection=1
\def\section#1{\bigbreak
 {\bf\number\numerodesection.\ \ #1}\nobreak\medskip
 \advance\numerodesection by1}

\mathcode`A="7041 \mathcode`B="7042 \mathcode`C="7043 \mathcode`D="7044
\mathcode`E="7045 \mathcode`F="7046 \mathcode`G="7047 \mathcode`H="7048
\mathcode`I="7049 \mathcode`J="704A \mathcode`K="704B \mathcode`L="704C
\mathcode`M="704D \mathcode`N="704E \mathcode`O="704F \mathcode`P="7050
\mathcode`Q="7051 \mathcode`R="7052 \mathcode`S="7053 \mathcode`T="7054
\mathcode`U="7055 \mathcode`V="7056 \mathcode`W="7057 \mathcode`X="7058
\mathcode`Y="7059 \mathcode`Z="705A


\def\hfl#1#2#3{\smash{\mathop{\hbox to#3{\rightarrowfill}}\limits
^{\textstyle#1}_{\textstyle#2}}}

\def\P{{\bf P}}

\def\Z{{\bf Z}}

\def\F{{\bf F}}

\def\GL{{\bf GL}}

\def\Hom{\mathop{\rm Hom}\nolimits}

\def\End{\mathop{\rm End}\nolimits}

\def\Ind{\mathop{\rm Ind}\nolimits}

\def\Card{\mathop{\rm Card}\nolimits}
\def\Gal{\mathop{\rm Gal}\nolimits}
\def\Ker{\mathop{\rm Ker}\nolimits}

\def\to{\rightarrow}
\def\lcm{\mathop{\rm lcm}\nolimits}

\def\mod{\mathop{\rm mod.}\nolimits}
\def\pmod#1{\;(\mod#1)}

\def\boxit#1{\vbox{\hrule\hbox{\vrule\kern1pt
       \vbox{\kern1pt#1\kern1pt}\kern1pt\vrule}\hrule}}
\def\cqfd{\hfill\boxit{\phantom{\i}}}

\newcount\numero \numero=1
\def\numeroter{{({\oldstyle\number\numero})}\ \advance\numero by1}

\def\pt{\hbox{\bf.}}

\newcount\refno 
\long\def\ref#1:#2<#3>{                                        
\global\advance\refno by1\par\noindent                              
\llap{[{\bf\number\refno}]\ }{#1} \pointir{\it #2} #3\goodbreak }

\newcount\refno 
\long\def\ref#1:#2<#3>{                                        
\global\advance\refno by1\par\noindent                              
\llap{[{\bf\number\refno}]\ }{#1} \pointir{\it #2} #3\goodbreak }

\def\citer#1(#2){[{\bf\number#1}\if#2\empty\relax\else,\ {#2}\fi]}

\def\pt{\hbox{\bf.}}

\def\boxit#1{\vbox{\hrule\hbox{\vrule\kern1pt
       \vbox{\kern1pt#1\kern1pt}\kern1pt\vrule}\hrule}}
\def\cqfd{\hfill\boxit{\phantom{\i}}}

\newbox\bibbox
\setbox\bibbox\vbox{\bigbreak
\centerline{{\pc BIBLIOGRAPHY}}
\nobreak
\ref{\pc BOURBAKI} (N):
Alg\`ebre, Chapitres 4 \`a 7,
<Masson, Paris, 1981, 422~p.>
\newcount\bourbaki \global\bourbaki=\refno

\ref{\pc DALAWAT} (C):
Serre's ``\thinspace formule de masse\thinspace'' in prime degree,
<Monats\-hefte Math.\ {\bf 166} (2012) 1, 73--92.
Cf.~arXiv\string:1004.2016v6.>
\newcount\monatshefte \global\monatshefte=\refno

\ref{\pc DALAWAT} (C) \& {\pc LEE} (JJ):
Tame ramification and group cohomology,
<J.\ Ramanujan Math.\ Soc.\ {\bf 32} (2017) 1, 51--74.
Cf.~arXiv\string:1305.2580v4.>
\newcount\dalawatlee \global\dalawatlee=\refno

\ref{\pc DALAWAT} (C):
Solvable primitive extensions,
<arXiv\string:1608.04673.>
\newcount\solprimp \global\solprimp=\refno

\ref{\pc DALAWAT} (C):
Wildly primitive extensions,
<arXiv\string:1608.04183.>
\newcount\wildprim \global\wildprim=\refno

\ref{\pc DEL \pc CORSO} (I) \& {\pc DVORNICICH} (R):
The compositum of wild extensions of local fields of prime degree,
<Monatsh.\ Math.\ {\bf 150} (2007) 4, 271--288.>
\newcount\delcorso \global\delcorso=\refno

\ref{\pc DEL \pc CORSO} (I), {\pc DVORNICICH} (R) \& {\pc MONGE} (M):
On wild extensions of a $p$-adic field,
<J.\ Number Theory {\bf 174} (2017), 322--342.  Cf.~aXiv\string:1601.05939.> 
\newcount\deldvomonge \global\deldvomonge=\refno

\ref{\pc DOERK} (K) \& {\pc HAWKES} (T):
Finite soluble groups,
<Walter de Gruyter \& Co., Berlin, 1992. xiv+891 pp.>
\newcount\doerkhawkes \global\doerkhawkes=\refno

\ref{\pc KOCH} (H):
Classification of the primitive representations of the Galois group of local fields,
<Invent.\ Math.\ {\bf 40} (1977) 2, 195--216.>
\newcount\koch \global\koch=\refno

\ref{\pc KOCH} (H):
On the local Langlands conjecture, <S\'eminaire de th\'eorie des
nombres de Grenoble {\bf 8} (1979-1980), 1--14.>
\newcount\kochllc \global\kochllc=\refno

\ref{\pc SERRE} (J-P):
Repr\'esentations lin\'eaires des groupes finis,
<Hermann, Paris, 1978, 182 p.>
\newcount\serrelinrep \global\serrelinrep=\refno

} 

\centerline{\bf $\F_p$-representations over $p$-fields}
\bigskip\bigskip 
\centerline{Chandan Singh Dalawat} 
\centerline{Harish-Chandra Research Institute}
\centerline{Chhatnag Road, Jhunsi, Allahabad 211019, India} 
\centerline{\tt dalawat@gmail.com}

\bigskip\bigskip

{{\bf Abstract}.  Let $p$ be a prime, $k$ a finite extension of $\F_p$
  of cardinal $q$, $l$ a finite extension of $k$ of group
  $\Sigma=\Gal(l|k)$, and $T$ a subgroup of $l^\times$.  Using the
  method of ``\thinspace little groups\thinspace'', we classify
  irreducible $\F_p$-representations of the group $G=T\times_q\Sigma$,
  the twisted product of $\Sigma$ with the $\Sigma$-module $T$.  We
  then use these results to classify irreducible continuous
  $\F_p$-representations of the profinite group $\Gal(\tilde K|K)$ of
  $K$-automorphisms of the maximal galoisian extension $\tilde K$ of a
  $p$-field $K$ with residue field~$k$.

\footnote{}{{\it MSC2010~:} Primary 11F80, 11S99, 20C20}
\footnote{}{{\it Keywords~:} Local fields, Galois representations,
  Little groups
}}

\bigskip
\bigbreak
{\bf 1.  Introduction}
\bigskip

\numeroter Let $p$ be a prime and let $K$ be a $p$-field, namely a
local field with finite residue field of characteristic~$p$.  Let
$\tilde K$ be a maximal galoisian extension of $K$.  Let $V$ be
the maximal tamely ramified extension of $K$ in $\tilde K$.  All
representations of the profinite groups $\Gal(\tilde K|K)$ and
$\Gal(V|K)$ appearing below are assumed to be {\it continuous\/}.  The
ramification group $\Gal(\tilde K|V)$, which is a pro-$p$-group, acts
trivially on any irreducible $\F_p$-representation of $\Gal(\tilde
K|K)$.  So classifying irreducible $\F_p$-representations of
$\Gal(\tilde K|K)$ comes down to classifying irreducible
$\F_p$-representations of $\Gal(V|K)$, which comes down to classifying
irreducible $\F_p$-representations of $\Gal(L|K)$ for every finite
tamely ramified galoisian extensions $L$ of $K$.

\numeroter For such $L$ with group of $K$-automorphisms
$G=\Gal(L|K)$ and inertia subgroup $G_0$, the projection $G\to G/G_0$
need not have a section, but $L$ has finite unramified extensions $L'$
for which the corresponding projections $G'\to G'/G_0'$ (where $G'$ is
$\Gal(L'|K)$ and $G_0'$ is its inertia subgroup) do have sections~;
the smallest such $L'$ is the one whose degree over $L$ is equal to
the order in $H^2(G/G_0,G_0)$ of the class of the extension $G$ of
$G/G_0$ by the $(G/G_0)$-module $G_0$~; see for example
\citer\dalawatlee(Lemma~2.3.4).  So it is enough to understand
irreducible $\F_p$-representations of $G$ in this {\it split\/} case~;
with a heavy heart, we choose a section of $G\to G/G_0$ in what
follows.

\numeroter Our treatment, which is completely
canonical and somewhat simpler than in the literature, is better
adapted to this arithmetic application because the inertia group $G_0$
does not come with a generator, only a canonical character
$\theta:G_0\to l^\times$, where $l$ is the residue field of $L$.  It
is based upon \S7 of \citer\koch(p.~205) and the method of
``\thinspace little groups\thinspace'' of Wigner and Mackey as exposed
in \S8.2 of \citer\serrelinrep(p.~62)~; I thank my friend UK
Anandavardhanan for pointing out the latter reference.  The material
is also worked out in \S4.1 of \citer\deldvomonge(p.~329).

\numeroter Let $k$ be the residue field of $K$~; the quotient
$\Sigma=G/G_0$ can be identified with $\Gal(l|k)$, which has the
canonical generator $\sigma:x\mapsto x^q$ ($x\in l$), where $q=\Card
k$.  The conjugation action of $\Sigma$ on $G_0$ is given by
$\sigma\pt t=t^q$ for every $t\in G_0$~; the character $\theta$ is
$\Sigma$-equivariant.  To determine the $\F_p$-representations of $G$,
we may forget the fields $K$ and $L$, and retain only $p$, $k$, $l$,
and $e=\Card G_0$.  This is done in \S2.  In \S3, we return to these
local fields and make an observation which will be useful
elsewhere \citer\wildprim().

\bigskip
\bigbreak
{\bf 2. Irreducible $\F_p$-representations of little groups}
\bigskip

\numeroter {\it Notation}.  Let us restart and rename.  Fix a prime
number $p$, fix a finite extension $k$ of $\F_p$, put $q=\Card k$, fix
a finite extension $l$ of $k$, put $f=[l:k]$, and denote by
$\sigma:x\mapsto x^q$ ($x\in l$) the canonical generator of
$\Sigma=\Gal(l|k)$.  Let $T\subset l^\times$ be a subgroup, let $e$ be
the order of $T$ (so that $q^f\equiv1\pmod e$), and let $\theta:T\to
l^\times$ be the inclusion (viewed as a character).  Finally, let
$G=T\times_q\Sigma$ be the twisted product of $\Sigma$ with the
$\Sigma$-module $T$ (for the galoisian action $\sigma\pt t=\sigma
t\sigma^{-1}=t^q$ ($t\in T$)).  Notice that the action is trivial, or
equivalently $G$ is commutative, if and only if $q\equiv1\pmod e$.

\numeroter {\it The problem}.  Classify irreducible
$\F_p$-representations of $G$.  

\numeroter {\it Notation}.  For every character $\chi:T\to l^\times$,
we denote by $d_\chi$ the order of $\chi$ and by $r_\chi$
(resp.~$s_\chi$) the order of the image $\bar p$ (resp.~$\bar q$) in
$(\Z/d_\chi\Z)^\times$.  Put $T_\chi=T/\Ker(\chi)$ and let
$\Sigma_\chi$ be the kernel of the action of $\Sigma$ on $T_\chi$~;
the group $\Sigma_\chi$ is generated by $\sigma^{s_\chi}$, and
$s_\chi$ is also the size of the $\Sigma$-orbit $\bar\chi$ of $\chi$.
We have the subgroup $G_\chi=T\times_q\Sigma_\chi$ of $G$ and the
quotient $\bar G_\chi=T_\chi\times\Sigma_\chi$ of $G_\chi$ which is
commutative by construction.  The numbers $d_\chi$, $r_\chi$, $s_\chi$,
the groups $\Sigma_\chi$, $G_\chi$, $\bar G_\chi$, and some characters
of these groups to be defined presently, depend only on the
$\Sigma$-orbit $\bar\chi$ of $\chi$~; we keep the notation light by
writing $\chi$ in the subscript instead of $\bar\chi$.  Let $\chi$
also stand for the faithful character $T_\chi\to l^\times$ coming from
$\chi$.

\numeroter We begin by determining the {\it irreducible $\tilde
  l$-representations of\/} $G$, where $\tilde l=l(\!\root{f'}\of1)$,
$f'=fp^{-v_p(f)}$, and $v_p(f)$ is the exponent of~$p$ in the prime
decomposition of~$f$~; they will turn out to be absolutely
irreducible.  These representations will be parametrised by pairs
$(\bar\chi,\lambda)$, where $\bar\chi\subset\Hom(T,l^\times)$ is the
$\Sigma$-orbit of a character $\chi:T\to l^\times$ (for the action
$\sigma:\chi\mapsto\chi^q$) and $\lambda\in\tilde l^\times$ is an
element of order dividing $fs_\chi^{-1}$.  This is achieved in several
steps.

\numeroter {\it The\/ $\tilde l$-representation\/
  $\rho_{\bar\chi,\lambda}$ of\/ $G$ associated to a pair
  $(\bar\chi,\lambda)$}.  Choose $\chi\in\bar\chi$, and let
  $\psi_{\chi,\lambda}:\Sigma_\chi\to\tilde l^\times$ be the unique character
  such that $\psi_{\chi,\lambda}(\sigma^{s_\chi})=\lambda$.  View the
  character $\chi\otimes\psi_{\chi,\lambda}$ of $\bar
  G_\chi=T_\chi\times\Sigma_\chi$ as a character of
  $G_\chi=T\times_q\Sigma_\chi$, and take the induced representation
  $\rho_{\bar\chi,\lambda}=\Ind^G_{G_\chi}(\chi\otimes\psi_{\chi,\lambda})$.

\numeroter To see that $\rho=\rho_{\bar\chi,\lambda}$ depends only on
the pair $(\bar\chi,\lambda)$, and for later use, let us make all this
explicit.  The quotient $G/G_\chi$ is generated by the image of
$\sigma$, so it can be identified with $\Z/s_\chi\Z$.  By definition,
the space $\Ind^G_{G_\chi}(\chi\otimes\psi_{\chi,\lambda})$ has an
$\tilde l$-basis $(b_i)_{i\in\Z/s_\chi\Z}$ on which the action of
$G_\chi$ is given by
$$
\rho( t)(b_i)=\chi(\sigma^i\pt t)b_i=\chi^{q^i}( t)b_i\qquad
(t\in T,\ i\in\Z/s_\chi\Z),
$$
and $\rho(\sigma^{s_\chi})(b_i)=\lambda b_i$.  This action is extended
to $G$ by $\rho(\sigma)(b_i)=b_{i+1}$ for
$i\not\equiv-1\pmod{s_\chi}$ and $\rho(\sigma)(b_{-1})=\lambda b_0$,
which gives back, as it should, the action of $\sigma^{s_\chi}$.  Now
it is clear that $\rho_{\bar\chi,\lambda}$ depends only on
$(\bar\chi,\lambda)$.

\numeroter {\it The\/ $\tilde l$-representation\/
  $\rho_{\bar\chi,\lambda}$ is absolutely irreducible and determines
  the pair\/ $(\bar\chi,\lambda)$}.  Write
  $\rho=\rho_{\bar\chi,\lambda}$.  None of the $T$-stable lines in
  $\rho$ (the $s_\chi$ lines on which $T$ acts respectively via the
  characters $\chi^{q^i}$, which are distinct for distinct
  $i\in\Z/s_\chi\Z$) is stable under $\sigma$ unless $s_\chi=1$, in
  which case $G_\chi=G$ and $\rho=\chi\otimes\psi_{\chi,\lambda}$, so
  $\rho$ is irreducible in every case, and in fact absolutely
  irreducible because the same agrument works over any finite
  extension of $\tilde l$.  Note that the $\Sigma$-orbit $\bar\chi$
  can be recovered from\/ $\rho$ because
  $\rho|_T=\oplus_{\eta\in\bar\chi}\;\eta$, and then $\lambda\in\tilde
  l^\times$ can be recovered because $\rho(\sigma^{s_\chi})$ is the
  homothety of ratio~$\lambda$.

\numeroter {\it Every irreducible\/ $\bar\F_p$-representations $\rho$ of\/ $G$ 
  come from a pair\/ $(\bar\chi,\lambda)$ as in\/ $\oldstyle(8)$}.
  Here $\bar\F_p$ is a maximal galoisian extension of $\tilde l$.  Let
  $\bar T=T/\Ker(\rho|_T)$, so that $\rho$ comes from an (irreducible)
  $\bar\F_p$-representation $\bar\rho$ of $\bar G=\bar
  T\times_q\Sigma$ whose restriction to $\bar T$ is faithful.  Let $P$
  be the intersection of the Sylow $p$-subgroups of $\bar G$.  The
  image $\bar\rho(P)$ is trivial because $P$ is a normal $p$-subgroup
  of $\bar G$, the characteristic of $\tilde l$ is $p$, and $\bar\rho$
  is irreducible \citer\serrelinrep(Chapitre~8, Proposition~26).  So
  $\bar\rho$ comes from a representation $\hat\rho$ of $\bar G/P$.
  Let $\Sigma'$ denote the kernel of the action of $\Sigma$ on $\bar
  T$, so that the subgroup $\bar G'=\bar T\times\Sigma'$ of $\bar G$
  is commutative.  As $\bar G'\cap P$ is the Sylow $p$-subgroup of
  $\Sigma'$, the order of $\hat\rho(\bar G')$ divides the order of
  $\tilde l^\times$, by our definition of $\tilde l$.  It follows that
  $\bar\rho|_{\bar G'}$ is a direct sum of characters $\bar
  G'\to\tilde l^\times$ such that the restriction to $\bar T$ of at
  least one of which --- call it $\xi$ --- is faithful

\numeroter View $\chi=\xi|_{\bar T}$ as a character of $T$.  Since the
order $e$ of $T$ divides the order $q^f-1$ of $l^\times$, we have
$\chi(T)\subset l^\times$.  Also, $\Ker(\chi)=\Ker(\rho|_T)$, so that
$\bar G'=\bar G_\chi=T_\chi\times\Sigma_\chi$, in the previous
notation.  Recall that $s_\chi$ is the size of the $\Sigma$-orbit
$\bar\chi$, and that $\Sigma_\chi$ is generated by $\sigma^{s_\chi}$.
Let $b_0\neq0$ be a vector (in the representation space of $\rho$) on
which $G_\chi$ acts through $\xi$, and define $\lambda\in\tilde
l^\times$ by $\xi(\sigma^{s_\chi})(b_0)=\lambda b_0$.  We claim that
$\rho=\rho_{\bar\chi,\lambda}$.

\numeroter  Put $b_i=\rho(\sigma^i)(b_0)$ for
$i\in[0,s_\chi[$.  Note that $\rho(\sigma)(b_{s_\chi-1})=\lambda b_0$
and 
$$
\rho( t)(b_i) =\rho( t\sigma^i)(b_0) =\rho(\sigma^i t^{q^i})(b_0)
 =\chi^{q^i}( t)b_i
 \qquad (t\in T,\,i\in[0,s_\chi[\,).
 $$
The characters $\chi^{q^i}$ are distinct for distinct $i\in[0,s_\chi[$,
therefore the family $(b_i)_{i\in[0,s_\chi[}$ is linearly independent.
Also, the subspace generated by the $b_i$ is $G$-stable, and in fact
equal to
$\rho_{\bar\chi,\lambda}=\Ind^G_{G_\chi}(\chi\otimes\psi_{\chi,\lambda})$
as described earlier.  Since $\rho$ is irreducible, we must have
$\rho=\rho_{\bar\chi,\lambda}$, as claimed.  Therefore~:


\numeroter {\it The set of irreducible\/ $\tilde l$-representations
  of\/ $G=T\times_q\Sigma$ is in natural bijection with the set of pairs\/
  $(\bar\chi,\lambda)$ consisting of the\/ $\Sigma$-orbit\/ $\bar\chi$
  of a character\/ $\chi:T\to l^\times$ and an element\/
  $\lambda\in\tilde l^\times$ of order dividing\/ $fs_\chi^{-1}$,
  where $s_\chi=\Card\bar\chi$. The pair\/ $(\bar\chi,\lambda)$ gives
  rise to the induced representation\/
  $\rho_{\bar\chi,\lambda}=\Ind_{G_\chi}^G(\chi\otimes\psi_{\chi,\lambda})$,
  where\/ $G_\chi=T\times_q\Sigma_\chi$, $\Sigma_\chi$ is generated by
  $\sigma^{s_\chi}$, and\/ $\psi_{\chi,\lambda}:\Sigma_\chi\to\tilde
  l^\times$ is the character such that\/
  $\psi_{\chi,\lambda}(\sigma^{s_\chi})=\lambda$.  All these
  representations are absolutely irreducible.} \cqfd

\numeroter {\it Some natural characters}.  The group $T$ comes
with the faithful character $\theta:T\to l^\times$, so $G$ has a
natural absolutely irreducible $\tilde l$-representation
$\rho_{\overline{\theta},1}$ of degree equal to the order of $\bar
q\in(\Z/e\Z)^\times$.  We could also consider $\theta^d$ for various
divisors $d$ of $e$.  Notice that $\theta$ allows us to identify
$\Hom(T,l^\times)$ with $\Z/e\Z$, and the $\Sigma$-orbit of
$\chi=\theta^i$ with the $\Sigma$-orbit of $i\in\Z/e\Z$ (for the
action $\sigma\mapsto(j\mapsto qj)$).

\smallbreak

\numeroter Let us now come to irreducible $\F_p$-representations $\pi$
of $G$, which are also treated in \citer\deldvomonge(Proposition~4.2).
The group $\Phi=\Gal(\tilde l\,|\F_p)$ acts on the set of irreducible
$\tilde l$-representations $\rho$ of $G$ by conjugation.  Let
$\varphi:x\mapsto x^p$ ($x\in\tilde l$) be the canonical generator of
$\Phi$.  If $\rho$ corresponds to the pair $(\bar\chi,\lambda)$ as
above, then ${}^\varphi\!\rho$ corresponds to the pair
$(\overline{\chi^p},\lambda^p)$.  The set of irreducible
$\F_p$-representations $\pi$ of $G$ is in natural bijection with the
set of $\Phi$-orbits $R$ for this action~; $\pi$ and $R$ correspond to
each other if $\pi\otimes_{\F_p}\!\tilde l=\oplus_{\rho\in R}\;\rho$.
If so, then $\deg\pi=(\deg\rho)(\Card R)$, for any $\rho\in R$.
Regarding this kind of ``\thinspace galoisian descent\thinspace'', see
for example \citer\bourbaki(V.60).

\numeroter Let us compute $\Card R$, or the size of the $\Phi$-orbit
$\overline{(\bar\chi,\lambda)}$ of any pair $(\bar\chi,\lambda)$ such
that $\rho_{\bar\chi,\lambda}\in R$.  Recall that $d_\chi$ is the
common order of every $\chi\in\bar\chi$, and that $r_\chi$
(resp.~$s_\chi$) is the order of $\bar p$ (resp.~$\bar q$) in
$(\Z/d_\chi\Z)^\times$.  (We already know that
$\deg\rho_{\bar\chi,\lambda}=s_\chi$).  The size of the $\Phi$-orbit
of the $\Sigma$-orbit $\bar\chi$ is $r_\chi s_\chi^{-1}$, so the
number of $\F_p$-conjugates of the pair $(\bar\chi,\lambda)$ (or the
size of the $\Phi$-orbit $R=\overline{(\bar\chi,\lambda)}$) is
$\lcm(r_\chi s_\chi^{-1},w_\lambda)$, where $w_\lambda$ is the degree
$[\F_p(\lambda):\F_p]$ (which obviously depends only on the
$\Phi$-orbit of $\lambda$), and the degree of $\pi$ is
$s_\chi\lcm(r_\chi s_\chi^{-1},w_\lambda)=\lcm(r_\chi,s_\chi
w_\lambda)$.  Therefore~:

\numeroter {\it The set of irreducible\/ $\F_p$-representations $\pi$
of\/ $G=T\times_q\Sigma$ is in natural bijection with the set of\/
$\Phi$-orbits\/ $R=\overline{(\bar\chi,\lambda)}$ of pairs\/
$(\bar\chi,\lambda)$ consisting of the\/ $\Sigma$-orbit\/ $\bar\chi$
of a character\/ $\chi:T\to l^\times$ and an element\/
$\lambda\in\tilde l^\times$ of order dividing\/ $fs_\chi^{-1}$,
where\/ $s_\chi$ is the order of\/ $\bar q\in(\Z/d_\chi\Z)^\times$
and\/ $d_\chi$ is the order of\/ $\chi$, under the correspondence\/
$\displaystyle\pi\otimes_{\F_p}\!\tilde
l=\oplus_{(\bar\chi,\lambda)\in R}\;\rho_{\bar\chi,\lambda}$, where\/
$\rho_{\bar\chi,\lambda}$ is the absolutely irreducible $\tilde
l$-representation of\/ $G$ attached to\/ $(\bar\chi,\lambda)$.  If\/
$r_\chi$ denotes the order of\/ $\bar p\in(\Z/d_\chi\Z)^\times$ and\/
$w_\lambda=[\F_p(\lambda):\F_p]$, then\/ $\deg\pi=\lcm(r_\chi,s_\chi
w_\lambda)$.}  \cqfd

\smallbreak 

\numeroter {\it The field of definition}.  Let
$k_{\bar\chi,\lambda}\subset\tilde l$ be the extension of $\F_p$ of
degree $\lcm(r_\chi s_\chi^{-1}, w_\lambda)$.  It follows for similar
reasons that there is a unique (absolutely) irreducible
$k_{\bar\chi,\lambda}$-representation $\rho'_{\bar\chi,\lambda}$ which
gives back $\rho_{\bar\chi,\lambda}$ upon changing the base to $\tilde
l$ in the sense that
$\rho'_{\bar\chi,\lambda}\otimes_{k_{\bar\chi,\lambda}}\tilde
l=\rho_{\bar\chi,\lambda}$~; we call $k_{\bar\chi,\lambda}$ {\it the
field of definition of\/} $\rho_{\bar\chi,\lambda}$ and henceforth
think of $\rho_{\bar\chi,\lambda}$ as a
$k_{\bar\chi,\lambda}$-representation.

\smallbreak

\numeroter We denote the irreducible $\F_p$-representation
of $G$ associated to the $\Phi$-orbit $\overline{(\bar\chi,\lambda)}$
by $\pi_{\overline{\bar\chi,\lambda}}$.  The degree of
$\pi_{\overline{\bar\theta,1}}$ is $r_\theta$ (the order of $\bar
p\in(\Z/e\Z)^\times$).  The notation is somewhat ambiguous because it
doesn't refer to the group $G$.  Indeed, if $l'$ is a finite extension
of $l$, then the same $\Phi$-orbit $\overline{(\bar\chi,\lambda)}$
also gives rise to an irreducible $\F_p$-representation
$\pi'=\pi_{\overline{\bar\chi,\lambda}}$ of the group
$G'=T\times_q\Sigma'$, where $\Sigma'=\Gal(l'|k)$.  The saving grace
is that if we use the galoisian projection $\gamma:\Sigma'\to\Sigma$
to view $G$ as a quotient of $G'$, then $\pi'=\pi\circ\gamma$.

\smallbreak

\numeroter We give some examples which will be useful later
\citer\wildprim() in classifying quartic extensions $E$ of a dyadic
field $K$ which have no intermediate quadratic extensions.  The set of
such $E$ was parametrised in
\citer\solprimp($\oldstyle14$) by the set of pairs $(\rho, D)$, where
$\rho$ is an irreducible $\F_2$-representation of $\Gal(\tilde K|K)$
(and $\tilde K$ is the maximal galoisian extension of $K$) and $D$ is
an $\F_2^2$-extension of the fixed field $F_\rho$ of the kernel of
$\rho$ such that $D$ is galoisian over $K$ and the resulting
conjugation action of $\Gal(F_\rho|K)$ on $\Gal(D|F_\rho)$ is given by
$\rho$.  If so, the group $\Gal(\hat E|K)$ (where $\hat E$ is the
galoisian closure of $E$ over $K$) is given by
$\F_2^2\times_\rho\Gal(F_\rho|K)$ \citer\solprimp().  Here we merely
construct all degree-$2$ irreducible $\F_2$-representations $\rho$ of
certain groups $G$ and identify the twisted product
$\F_2^2\times_\rho G$.  Why it suffices to consider only these $G$ was
explained in \citer\solprimp($\oldstyle15$)
and will also become clear at the very end.

\numeroter {\it A\/ $(\Z/3\Z)$- and an\/ ${\goth A}_3$-example}.    
Take
$p=2$, $f=3$, $e=1$, so that $T=\{1\}$, and $G=\Sigma=\Z/3\Z$.  We
have $\tilde l=l(\!\root3\of1)$.  The only $\chi:T\to l^\times$ is the
trivial character~$1$~; for it, there are three possible $\lambda$,
namely $1$, $\root3\of1$ and ${\root3\of1}^2$.  So we get three
$\tilde l$-characters, namely $\rho_{\bar1,1}$,
$\rho_{\bar1,\root3\of1}$ and $\rho_{\bar1,{\root3\of1}^2}$ of which
the latter two are in the same $\Phi$-orbit, and these are the only
irreducible $\tilde l$-representations of $G$.  Thus we get two
$\F_2$-representations, namely the trivial representation
$\pi_{\overline{\bar1,1}}$ and the irreducible degree-$2$
representation $\pi=\pi_{\overline{\bar1,\root3\of1}}$~; the latter is
not absolutely irreducible.  The group $\F_2^2\times_\pi G$ is
isomorphic to ${\goth A}_4$.

Or keep $p=2$ and take $q\equiv1\pmod3$, $f=1$, $e=3$, so that $\tilde
l=l$.  Then $G=T={\goth A}_3$ is cyclic of order~$3$, the three
irreducible $l$-representations are $\rho_{\bar1,1}$,
$\rho_{\bar\theta,1}$, $\rho_{\bar\theta^2,1}$ (all three of degree~1)
of which the latter two are in the same $\Phi$-orbit, so the two
irreducible $\F_2$-representations are $\pi_{\overline{\bar1,1}}$
(trivial) and $\pi=\pi_{\overline{\bar\theta,1}}$ (degree~2).  The
group $\F_2^2\times_\pi G$ is isomorphic to ${\goth A}_4$, as before.

\smallbreak

\numeroter {\it An\/ ${\goth S}_3$-example}.  Keep $p=2$
and take $q\equiv-1\pmod3$, $f=2$, $e=3$, so that $\Sigma=\Z^\times$,
$G$ is isomorphic to ${\goth S}_3$, and $\tilde l=l$.  The only
characters $T\to l^\times$ are $1$ (of order $d=1$) and $\theta$,
$\theta^2$ (of order $d=3$)~; they fall into two $\Sigma$-orbits,
namely $\bar1$ (of size $s=1$) and $\bar\theta$ (of size $s=2$).  The
only possible $\lambda$ in either case is $\lambda=1$.  So
$\rho_{\bar1,1}$ (of degree~$1$) and $\rho_{\bar\theta,1}$ (of
degree~$2$) are the only two (absolutely) irreducible
$l$-representations of $G$, each of which is its own $\Phi$-orbit.
Therefore there are two irreducible $\F_2$-representations of $G$,
namely $\pi_{\overline{\bar1,1}}$ (the trivial representation) and
$\pi=\pi_{\overline{\bar\theta,1}}$ (of degree~$2$).  In fact,
$\pi:G\to\GL_2(\F_2)$ is an isomorphism, and $\F_2^2\times_\pi G$ is
isomorphic to ${\goth S}_4$.

\numeroter Another general algebraic observation we need is the
following lemma culled from \citer\deldvomonge(4.9)~; see also
\citer\doerkhawkes(p.~154).  Let $G$ be any finite group, $F$ any field, $E$
a finite galoisian extension of $F$, $W$ an {\it absolutely\/}
irreducible $E$-representation of $G$ such that the conjugates
${}^\sigma W$ (\hbox{$\sigma\in\Gal(E|F)$}) of $W$ are all
inequivalent.  By galoisian descent, there is a unique (irreducible)
$F$-representation $V$ of $G$ such that $V\otimes_F
E=\oplus_{\sigma\in\Gal(E|F)}{}^\sigma W$.  By Schur's lemma, we have
$\End_{E[G]}(W)=E$ and also $\End_{F[G]}(V)=E$.

Let $m>0$ be an integer.  For every $a=(a_i)_{i\in[1,m]}$ in $E^m$, we
have the $F[G]$-morphism $\varphi_a:V\to V^m$ sending $x$ to
$(a_ix)_{i\in[1,m]}$~; it is injective if and only if $a\neq0$.  For
$a\neq0$, the image $\varphi_a(V)$ depends only on the line $\bar
a\subset E^m$ generated by $a$.

\numeroter {\it The map\/ $\bar a\mapsto\varphi_a(V)$ is a bijection
of the set\/ $\P_{m-1}(E)$ of lines in\/ $E^m$ with the set of
submodules of\/ $V^m$ isomorphic to\/ $V$.  In particular, if\/ $E$ is
finite and\/ $q=\Card(E)$, then the number of such submodules is\/
$(q^m-1)(q-1)^{-1}$.}

{\it Proof}.  The map in question is injective~: indeed, if $a\neq0$
and $b\neq0$ are in $E^m$, and if $\varphi_a(V)=\varphi_b(V)$, then
(slightly abusing notation) $\varphi_b^{-1}\circ\varphi_a$ is a
$G$-automorphism of $V$, so a homothety of some ratio $\xi\in
E^\times$, therefore $a=\xi b$ and $\bar a=\bar b$.  Next, the map
$\bar a\mapsto\varphi_a(V)$ is surective~: if $\psi:V\to V^m$ is an
injective $G$-morphism, the maps $\pi_i\circ\psi$, where the
$\pi_i:V^m\to V$ are the canonical projections, are homotheties of
some ratio $a_i\in E$ such that $a=(a_i)_{i\in[1,m]}$ is $\neq0$,
and $\psi=\varphi_a$.  \cqfd

\bigskip
\bigbreak
{\bf 3.  Irreducible $\F_p$-representations over $p$-fields}
\bigskip

\numeroter Let $K$ be a $p$-field, $k$ its residue
field, $q=\Card k$, and let $V_0$ (resp.~$V$, resp.~$\tilde K$) be the
maximal unramified (resp.~tamely ramified, resp.~galoisian) extension
of $K$. Put $\Gamma_0=\Gal(V_0|K)$ and $\Gamma=\Gal(V|K)$.  We have
seen in the Introduction that every irreducible $\F_p$-representation
of $\Gal(\tilde K|K)$ factors through $\Gamma$.

For every $n>0$, put $e_n=p^n-1$ and $V_n=V_0(\!\root e_n\of\varpi)$,
where $\varpi$ is a uniformiser of $K$.  It doesn't matter which
$\varpi$ we choose because $V_n$ is also obtained by adjoining the
family $\root e_n\of x$ (indexed by $x\in V_0^\times$) to $V_0$ .
Every $V_n$ is galoisian over $K$~; put $\Gamma_n=\Gal(V_n|K)$, so
that $V=\lim\limits_{\longrightarrow}V_n$ and
$\Gamma=\lim\limits_{\longleftarrow}\Gamma_n$.  The salient quotients
$\Gamma_n$ of $\Gamma$ have nothing to do with the ramification
filtration on $\Gamma$ which is quite simply $\Gamma^0\subset\Gamma$,
where $\Gamma^0=\Gal(V|V_0)$ is the inertia subgroup.  If $p=2$, then
$V_1=V_0$.  Note that if $K$ has characteristic~$0$, then the
$p$-torsion subgroup ${}_pV_n^\times$ of $V_n^\times$ has order~$p$
(because $V_1$ contains $\root{p-1}\of{-p}\,$).


\numeroter Note that $V_n$ is the compositum of all finite
extensions of $K$ of ramification index dividing $e_n$, so the
indexing has something to do with ramification afterall.  Note also
that if $a=v_p(q)$ is the exponent of $p$ in $q$, then $V_a$ is the
maximal tamely ramified abelian extension of $K$.
 
\numeroter {\it Let\/ $n>0$ be an integer.  Every irreducible\/
$\F_p$-representation of\/ $\Gamma$ of degree dividing\/~$n$ factors
through the quotient\/ $\Gamma_n$ of\/ $\Gamma$.}

{\it Proof}. Let $\pi$ be such a representation, and let $L$ be a
finite unramified extension of the fixed field $V^{\Ker(\pi)}$ which
is split over $K$ in the sense the inertia subgroup $G_0$ of
$G=\Gal(L|K)$ has a complement in $G$~; by hypothesis, $\pi|_{G_0}$ is
faithful.  It suffices ($\oldstyle27$) to show that the ramification
index~$e$ of $L$ over $K$ divides $e_n$.

Let $l$ be the residue field of $L$.  The filtration $G_0\subset G$ is
split by hypothesis~; the choice of a section $G/G_0\to G$ leads to an
isomorphism of $G$ with $T\times_q\Sigma$, where $\Sigma=\Gal(l|k)$
and $T\subset l^{\times}$ is the subgroup of order $e$.  Since
$\pi|_T$ is faithful, $\chi|_T$ is faithful for any character
$\chi:T\to l^\times$ which occurs in $(\pi|_T)\otimes\tilde l$ as in
($\oldstyle11$).  Therefore the order of $\chi$ is $e$.  Let $r$ be
the order of $\bar p\in(\Z/e\Z)^\times$, so that $p^r\equiv1\pmod e$.
Since $n$ is a multiple of $r$ ($\oldstyle19$), we have
$p^n\equiv1\pmod e$, and hence $e$ divides $e_n=p^n-1$.  \cqfd

\numeroter Since there are only finitely many irreducible
$\F_p$-representations $\pi$ of $\Gamma$ of given degree~$n$ (because
$\Gamma$ is finitely generated and $\GL_n(\F_p)$ is finite), there are
finite extensions $M$ of $K$ such that every irreducible
$\F_p$-representation of $\Gamma$ of degree~$n$ factors through
$\Gal(M|K)$.  For $n=1$, the smallest possible $M$ is clearly
$L_1=K(\!\root{p-1}\of{K^\times})$, which was used in
\citer\monatshefte() (and in \citer\delcorso() in the
characteristic-$0$ case)~; recently I've discovered that this
observation was made already in \citer\kochllc(V.9).  

\numeroter {\it A partition}.  Notice that {\it ramified\/} irreducible
 $\F_p$-representations $\pi$ of $\Gamma$ of degree $n>0$ can be
 partitioned into classes labelled by the divisors $r$ of $n$.  The
 representation $\pi$ belongs to the class labelled by $r$ if $r$ is
 the smallest (in the sense of divisibility) divisor of $n$ such that
 $\pi$ factors through $\Gamma_r$~; equivalently, $r$ is the order of
 $\bar p\in(\Z/e\Z)^\times$, where $e$ is the ramification index over
 $K$ of the fixed field $V^{\Ker(\pi)}$.  If $p=2$, then the class of
 label~$1$ is~$\hbox{\O}$ because $\Gamma_1=\Gamma_0$ and $\pi$ would
 be unramified.
 

\numeroter For every $n>0$, $K_n=K(\!\root e_n\of1)$ is the unramified
 extension of $K$ of degree equal to the order $s_n$ of $\bar
 q\in(\Z/e_n\Z)^\times$~; put $L_n=K_n(\!\root e_n\of{K_n^\times})$,
 so that $L_n\subset V_n$ and indeed $V_n=L_nV_0$.  Note that $L_n$ is
 the maximal abelian extension of $K_n$ of exponent dividing~$e_n$, so
 it is galoisian over $K$~; put $G_n=\Gal(L_n|K)$.  We have
$$
V_0=\lim\limits_{\longrightarrow}K_n,\quad
V=\lim\limits_{\longrightarrow}L_n,\quad
\Gamma=\lim\limits_{\longleftarrow} G_n.
$$
Note that if $K$ has characteristic~$0$, then ${}_pL_n^\times$ has
order~$p$ (because $L_1$ contains $\root{p-1}\of{-p}\,$).

In the proof of the next proposition, we shall need to consider
 certain finite galoisian extensions $M$ of $K$.  We denote the
 residue fields of $K_n$ (resp.~$L_n$, resp.~$M$) by $k_n$
 (resp.~$l_n$, resp.~$m$).  Note that $L_n$ is split over $K$ because
 $L_n=L_{n,0}(\!\root e_n\of\varpi)$ for any uniformiser $\varpi$ of
 $K$, where $L_{n,0}$ is the maximal unramified extension of $K$ in
 $L_n$.  If $M$ is unramified over $L_n$, then $M$ is also split over
 $K$.

\numeroter {\it   Every irreducible\/
$\F_p$-representation of\/ $\Gamma$ of degree dividing\/~$n$ factors
through the finite quotient\/ $G_n=\Gal(L_n|K)$ of\/ $\Gamma$.}

{\it Proof}.  Let $\pi'$ be such a representation, and recall that it
factors through $\Gamma_n$ ($\oldstyle28$).  Let $M$ be a finite
unramified extension of $L_n$ such that $H=\Gal(M|K)$ has the property
claimed for $G_n$~; we have to show that
$\Gal(M|L_n)\subset\Ker(\pi')$.  Choose a uniformiser $\varpi$ of $K$
and identify $H$ with $T\times_q\Gal(m|k)$ and $G_n$ with
$T\times_q\Gal(l_n|k)$ as above~; these identifications are compatible
with the galoisian projections $\gamma:H\to G_n$ and
$\Gal(m|k)\to\Gal(l_n|k)$.  The representation
$\pi'=\pi_{\overline{\bar{\eta},\mu}}$ of $H$ in question is
associated to the $\Phi$-orbit of a pair $(\bar{\eta},\mu)$ consisting
of the $\Gal(m|k)$-orbit of some character $\eta:T\to m^\times$ and
some $\mu\in\tilde m^\times$ as in ($\oldstyle19$).  We certainly have
$\eta(T)\subset k_n^\times$, therefore $\bar{\eta}$ can be viewed as
the $\Gal(l_n|k)$-orbit $\bar\chi$ of a character $\chi:T\to
l_n^\times$.  Recall that the degree $w_{\mu}=[\F_p(\mu):\F_p]$
divides~$n$, therefore the order of $\mu$ divides $e_n$, and hence
$\mu\in k_n^\times$~; call it $\lambda$ in this avatar.  The
$\Phi$-orbit of this new pair $(\bar\chi,\lambda)$ gives a
representation $\pi=\pi_{\overline{\bar\chi,\lambda}}$ of $G_n$, and
$\pi'$ factors through it ($\pi'=\pi\circ\gamma$), as in
($\oldstyle21$).  \cqfd

\numeroter {\it Remark}.  We define the {\it optimal\/} quotient of
$\Gamma$ in degree~$n$ to be the smallest quotient of $\Gamma$ through
which every irreducible $\F_p$-representation $\pi$ of $\Gamma$ of
degree~$n$ factors.  We also say that the corresponding extension $M_n$
of $K$ is the {\it optimal\/} extension in degree~$n$~; it is the
compositum, over all $\pi$, of the extensions $V^{\Ker(\pi)}$ of $K$.
For $n=1$, the quotient $G_1$ of $\Gamma$ and the corresponding
extension $L_1=K(\!\root p-1\of{K^\times})$ of $K$ are clearly
optimal. For $n>1$, the extension $M_n$ is introduced
in \citer\deldvomonge() but we prefer working with $L_n$ (which
contains $M_n$ by ($\oldstyle33$)) because $L_n$ is very explicit,
contains $\root p\of1$ in characteristic~$0$, and it is split over $K$
so that the general theory of little groups as explained above can be
applied.  Note that $M_n$ contains the unramified extension of $K$ of
degree~$e_n$ and every totally ramified extension of $K$ of degree
dividing~$e_n$.  We haven't checked whether $M_n=L_n$ in general.

\numeroter {\it The\/ $(\Z/3\Z)$-, ${\goth A}_3$- and\/ ${\goth
S}_3$-examples}.  Consider the case $p=2$ and $n=2$, so that $e_2=3$.
If $q\equiv1\pmod3$, then we have $K_2=K$ and
$L_2=K(\root3\of{K^\times})$, so that $L_2$ contains the unramified
cubic extension and the three ramified cubic extensions (all three
cyclic) of $K$.  We thus get back the $(\Z/3\Z)$- and ${\goth
A}_3$-examples of ($\oldstyle23$).  If $q\equiv-1\pmod3$, then
$[K_2:K]=2$, so that $L_2$ contains the unramified cubic extension and
the unique ${\goth S}_3$-extension (namely
$K(\root3\of1,\root3\of\varpi)$, where $\varpi$ is a uniformiser) of
$K$.  We thus get back the $(\Z/3\Z)$-example of ($\oldstyle23$) and
the ${\goth S}_3$-example of ($\oldstyle24$).


\bigbreak\bigbreak
\unvbox\bibbox 

\bye